\documentclass[a4paper,12pt]{article}
\usepackage{citesort,geometry,enumerate,amsmath,amssymb,latexsym,theorem,calc}
{\theorembodyfont{\rmfamily}\newtheorem{example}{Example}[section]}
{\theorembodyfont{\rmfamily}\newtheorem{defn}[example]{Definition}}
\newtheorem{prop}[example]{Proposition}

\newtheorem{propdef}[example]{Proposition and Definition }

{\theorembodyfont{\rmfamily}\newtheorem{rem}[example]{Remark}}
\newtheorem{cor}[example]{Corollary}

\newcommand{\Reno}{{\mathbb R}}

\newenvironment{pf}{{\bf Proof:}}{\hfill $\Box$

\mbox{}}
 \newcommand{\st}[1]{\mathrm{St}_{#1}\,}

\def\leq{\leqslant}

\def\le{\leqslant}
\def\C{\mathcal{C}}
\def\ob{\mathsf{Ob}}
\def\io{^{-1}}
\def\cE{\mathcal{E}}
\def\Sym{\mathsf{Sym}}
\def\L{\mathcal{L}}

\def\U{\mathcal{U}}
\def\subs{\subseteq}
\def\St{\mathrm{St}}
\def\Loc{\mathbf{Loc}}

\begin{document}

\title{\Large \bf Local subgroupoids II: Examples and properties
 }    \author{
 Ronald Brown  \\ School of Informatics \\ Mathematics
Division
\\ University of Wales  \\ Bangor, Gwynedd \\ LL57 1UT, U.K.
\\r.brown@bangor.ac.uk \\ \and
\.{I}lhan \.{I}\c{c}en    \\  University of  \.{I}n\"{o}n\"{u} \\
Faculty of Science and Art
\\ Department of Mathematics
\\ Malatya/ Turkey \\ iicen@inonu.edu.tr
\\ \and Osman Mucuk \\ University of Erciyes\\
Faculty of Science and Art \\Department of
Mathematics\\Kayseri/Turkey\\ mucuk@erciyes.edu.tr} \maketitle

\begin{center}{\bf
UWB Maths Preprint 00.15}
\end{center}
\begin{abstract}
The notion of local subgroupoid as a generalisation of a local
equivalence relation was defined in a previous paper by the first
two authors. Here we the notion of star path connectivity for a
Lie groupoid to  give an important new class of examples,
generalising the local equivalence relation of a foliation, and
develop in this new context basic properties of coherence, due
earlier to Rosenthal in the special case. These results are
required for further applications to holonomy and monodromy.
\end{abstract}

\section*{Introduction}

Any foliation gives rise to a local equivalence relation, defined
by the path components of local intersections of small open sets
with the leaves. Local equivalence relations were generalised to
local subgroupoids in  a previous paper by the first two authors
\cite{Br-Ic}, referred to hereafter as paper I.  In this paper we
show that a basic topological groupoid notion, that of {\em
identity star path component}, leads easily to a local subgroupoid
of a wide class of Lie groupoids and this generalises the local
equivalence relation of a foliation. We define local subgroupoids
$c_1(Q,\U)$ for certain open covers $\U$ of the object space of a
Lie groupoid $Q$. Further, we show that the theory of coherence,
which is prominent in the papers of Rosenthal \cite{Ro1,Ro2},
generalises nicely to the local subgroupoid case. The application
of  paper I to the holonomy and monodromy Lie groupoids of local
subgroupoids (to be dealt with elsewhere) requires a condition on
the local subgroupoid of having a `globally adapted atlas'. We
develop the `coherence' theory for giving conditions on $\U$ for
it to be globally adapted to the local subgroupoid $c_1(Q,\U)$
(Corollary \ref{glo-adapt}). This is related to the construction
of a locally Lie groupoid from a foliation \cite{Ku1,Br-Mu2}.

We also develop similar theory for Lie groupoids with a `path
connection'  \cite{Br-S,Ma,Vi} $\Gamma$ leading to a local
subgroupoid $c_{_\Gamma}(Q,\U)$ for certain open covers $\U$.

This paper is strongly influenced by papers of Rosenthal
\cite{Ro1,Ro2} on local equivalence relations, a concept  due
originally to Grothendieck and Verdier \cite{Gr-Ve} in a series of
exercises presented as open problems concerning the construction
of a certain kind of topos. The concept was investigated more
recently by Kock and Moerdijk \cite{Ko-Mo1,Ko-Mo2}.  The main aims
of the papers \cite{Gr-Ve,Ko-Mo1,Ko-Mo2,Ro1,Ro2} are towards the
connections with sheaf theory and topos theory.

The starting point of our investigation was to notice that an
equivalence relation on $X$ is a wide subgroupoid of the groupoid
$X \times X$. However an equivalence relation is just one of the
standard examples of a groupoid, and so it natural to consider the
corresponding theory for subgroupoids of a given groupoid $Q$. The
expectation is that this will allow applications to combinations
of foliation and bundle theory,  since a standard example of a Lie
groupoid is the Ehresmann symmetry groupoid of a principal bundle
\cite{Ma}.

In the case $Q$ is the indiscrete groupoid $X \times X$, we
recover the well known concept of {\em local equivalence
relation},

The remarkable fact is that the theory goes very smoothly, and so
suggests it is a natural generalisation of the foliation case, and
one which illuminates some constructions in that area.

%------------------------------------------------------------------------------
\section{Local subgroupoids: definitions and examples}

We first recall some definitions from \cite{Br-Ic}.

Consider a groupoid $Q$ on a set $X$ of objects, and suppose also
$X$ has a topology. For any open subset $U$ of $ X$ we write $Q|U$
for the full subgroupoid of $Q$ on the object set $U$. Let
$L_Q(U)$ denote the set of all subgroupoids of $ Q|U$ with object
set $U$ (these are called {\em wide} subgroupoids of $Q|U$). For
$V\subseteq U$, there is a restriction map $L_{UV}\colon L_Q(U)
\to L_Q(V)$ sending $H$ in $L_Q(U)$ to $H|V$. This gives $L_Q$
the structure of presheaf on $X$.

We interpret  the sheaf $p_Q:\L _Q \to X$ constructed in the
usual way  from the presheaf $L_Q$.

For $x\in X$, the stalk ${p_Q}^{-1}(x)$ of $\L _Q$ has elements
the germs $[U, H_U]_x$ where $U$ is open in $X$, $x\in U$, $H_U$
is a wide subgroupoid of $Q|U$, and the equivalence relation
$\sim_x$ yielding the germs at $x$ is that $H_U\sim_x K_V$, where
$K_V$ is wide subgroupoid of $Q|V$, if and only if there is a
neighbourhood $W$ of $x$ such that $W\subseteq U\cap V$ and $H_U|W
=K_V|W$. The topology on $\L_Q$ is the usual sheaf topology, with
a sub-base of sets $\{ [U,H]_x:x \in U \}$ for all open $U$ of $X$
and wide subgroupoids $H$ of $G|U$.

\begin{defn}
A {\it local subgroupoid} of $Q$ on the topological space $X$ is a
continuous global section of the sheaf $p_Q:\L _Q\to X$ associated
to the presheaf $L_Q$.
\end{defn}

Two standard  examples of $Q$ are   $Q=X$, $Q=X\times X$, where
$X\times X$ has the multiplication $(x,y)(y,z)=(x,z)$.  In the
first case, $L_X$ is already a sheaf and $\L _X\to X$ is a
bijection. More generally, we have:
\begin{prop} If $Q$ is a bundle of groups, then $L_Q$ is a sheaf.
\end{prop}
\begin{pf}
By our assumption,  if $U$ is a subset of $X$ then a wide
subgroupoid $H|U$ of $Q|U$ is uniquely defined by the values
$H(x)$ for all $x \in U$. This easily implies the usual two
compatibility conditions for a sheaf.
\end{pf}
In the case $Q$ is the indiscrete groupoid $X\times X$ the local
subgroupoids of $Q$ are the local equivalence relations on $X$,
as mentioned in the Introduction. It is known that $L_{X\times
X}$ is in general not a sheaf \cite{Ro1}.
\begin{defn}\label{loc}
If $G$ is a wide subgroupoid of the groupoid $Q$ on $X$, then  $loc(G)$ is the
local subgroupoid defined by
\[             loc(G)(x) =[X, G]_x.    \]
\end{defn}

This gives a wide and important class of local subgroupoids, but
we are more interested in those which derive from connectivity
considerations on a topological groupoid. For this we need to
discuss atlases for local subgroupoids, which is the major way of
giving a local subgroupoid.

\begin{defn}An {\it atlas } ${\cal U}_H = \{(U_i, H_i):i\in I\}$ for a
local subgroupoid  consists of an open cover ${\cal U}= \{U_i:i\in
I\}$ of $X$, and for each $i\in I$ a wide subgroupoid $H_i$ of
$Q|U_i$  such that the following compatibility condition holds:
\begin{enumerate}
  \item []{\bf
Comp$(H)$:} for all $i,j \in I,x\in U_i\cap U_j$  there is an
open set $W$ such that $x \in W \subs U_i\cap U_j $ and $H_i|W =
H_j|W$.
\end{enumerate}
The {\em local subgroupoid $s$ of the atlas} is then well defined
by $s(x)=[U_i, H_i]_x, \, x \in X$.

The above atlas is {\em compatible} with an atlas  ${\cal U'}_{H'}
= \{(U'_j, H'_j):j\in J\}$ if for all $i \in I, j \in J$ and $x
\in U_i \cap U'_j$ there is an open set $W$ such that $x \in W
\subs U_i\cap U'_j $ and $H_i|W = H'_j|W$. Clearly, two compatible
atlases define the same local subgroupoid.
\end{defn}
It is well known from general sheaf theory that any local
subgroupoid has a compatible atlas. Note also that the atlas
$\{(X,H)\}$ with a single element determines the local subgroupoid
$loc(H)$. So the atlas is a crucial part of the construction of a
local subgroupoid  $s$.

\section{The star path component of a topological groupoid}

A key concept for topological groups is the path component of the
identity. The analogue for topological groupoids is the star
identity path component.

\begin{defn}Let $Q$ be a topological groupoid. If $x \in \ob(Q)$ we write
$\st{Q}x$ for the star of $Q$ at $x$, namely the union of all the
$Q(x,y)$ for $ y \in \ob (Q)$. The {\it star identity path
component} $C_1(Q)$ of $Q$ consists of all $ g \in Q$ such that if
$x = \alpha (g)$ then there is a path in $\st{Q}x$ joining $g$ to
the identity $1_x$. Such a path is called a {\em star path}. We
say $Q$ is {\em star path connected } if $Q= C_1(Q)$.
\end{defn}

\begin{prop}
The star identity path component of $Q$ is a  subgroupoid of $Q$.
\end{prop}
\begin{pf}
Write $C$ for this star identity component. Let $g\in Q(x,y), \; h
\in Q(y,z)$ and suppose also $g,h \in C$. Then there are paths
$l_t$ in $\st{Q}x$, $m_t$ in $\st{Q}y$ such that $l_0=g, l_1=1_x,
m_0=h,m_1=1_y$. Hence $g.m_t$ is a path in $\st{Q}x $ joining $gh$
to $g$, and this composed with $l_t$ joins $gh$ to $1_x$. So $C$
is closed under composition.

If $ g \in C(x,y)$ and $l_t$ joins $g$ to $1_x$ then $g\io
l_{1-t}$ joins $g \io$ to $1_y$. So $C$ is a subgroupoid of $Q$.
\rule{1cm}{0cm}\end{pf}

Note that \cite[Example II.3.7, p46]{Ma} gives an example where
$C_1(Q)$ is not normal in $Q$.

We will later need the following.

\begin{prop} \label{generate} Let $Q$ be star path connected and let $\U$ be an
open cover of $X$. Then $Q$ is generated by the subgroupoids
$C_1(Q|U)$ for all $U \in \U$.
\end{prop}
\begin{pf}
Let $g \in \St _Qx$. Then there is a path $\lambda$ in $\St _Qx$
from $g$ to $1_x$. Let $\mu=\beta \lambda$. By the Lebesgue
covering lemma, we can write $\mu=\mu_1 + \cdots +\mu_n$ where
each $\mu_r$ lie in an open set $U_r$ of $\U$ for $r=1,
\ldots,n$. Then we can write $\lambda=\lambda_1 + \cdots
+\lambda_n$ where $\beta \lambda_r=\mu_r,r=1, \ldots,n$, and
$\lambda_r$ is a path in $\St_Qx \cap U_r$ for $r=1, \ldots,n$.
Let $g_r= -\lambda_r(0)+ \lambda_r(1)$. Then $g_r \in C_1(Q|U_r)$
and $g=g_1+\cdots+g_n$.
\end{pf}

Now  we recall some major examples of Lie groupoids.

\begin{example} Let $\cE$ be a principal bundle $p: E \to B$
with group $\Omega$. Then $E \times E$ is certainly a topological
groupoid, and so also is its quotient $Q=E \times_G E$ by the
diagonal action of the topological group $\Omega$. If $b,b' \in
B$, we can by choosing a point in $p\io(b)$, identify $Q(b,b')$
with the $\Omega$-maps $p\io(b) \to p\io(b')$. For this reason, we
also write $\Sym(\cE)$ for $Q$. In the case $\cE$ is locally
trivial, and assuming $X$ is a manifold, the topology may also be
constructed from this alternative description, since the
$\Omega$-maps $p\io(b) \to p\io(b')$ may, again by choosing a
point in $p\io(b)$, be identified with the elements of $\Omega$.
It is this description we now use.

Consider in particular the double cover of the circle $p: S^1 \to
S^1$ given by $z \mapsto z^2$. In this case $\Omega$ is the cyclic
group of order 2.  Let $Q=\Sym(p)$.

This groupoid $Q$ is star path connected. For suppose $g \in
Q(z,w)$. Let $\lambda$ be a path of shortest length in $S^1$ from
$z$ to $w$ (if $z=-w$ then there are two such paths). Let $u \in
S^1$ satisfy $u^2=z$. Since $p$ is a covering map, there are
unique paths $\lambda^+,\lambda^-$ starting at $u,-u$ and covering
$\lambda$. Let $v=\lambda^+(1)$. Then $g$ is a bijection $\{u,-u\}
\to \{v,-v\}$. If $g(u)=v$, then the pair of paths $\lambda^+,
\lambda^-$ define a star path from  the identity on  $\{u,-u\}$ to
$g$. If  $g(u)=-v$, then such a path is determined by $\lambda'$,
the shortest path joining $z$ to $w$ in the opposite direction
round $S^1$, and its corresponding lifts.

However, if $U=S^1$ with a single point removed, then $Q|U$ is not
star-connected, since if $\lambda$ is a path joining $z$ to $w$ in
$U$, then $\lambda'$ is not a path in $U$. \hfill $\Box$
\label{sym}
\end{example}

\begin{example}
Let $\Omega$ be a Lie group acting smoothly on the right of a
$\C^r$-manifold $X$. Form the Lie action groupoid $Q=X \rtimes
\Omega$. Even if $Q$ is star path connected, this is not
necessarily so for $Q|U$ for all open subsets $U$ of $X$.\hfill
$\Box$
\end{example}

\section{Local subgroupoids and star path connectivity}

The previous notions give us our major examples of new and
interesting local subgroupoids.

\begin{example}Consider  an equivalence relation $E$ on the space $X$. Then
for each open set $U$ of $X$ we have an equivalence relation $E|U$
on $U$ and  we can consider the partition of $U$ given by the path
components of the equivalence classes of $E|U$. In general, this
will not give us a local equivalence relation. Instead we need to
assume given  an open cover ${\cal U}= \{U_i: i \in I \}$ of $X$
satisfying the compatibility condition that for all $i,j \in
I,x\in  U_i\cap U_j$  there is an open set $W$ such that $x \in W
\subs U_i\cap U_j $ and the path components of $E|W$ are the
intersections with $W$ of the path components of the classes of
each of $E|U_i, \;E|U_j$. The resulting local equivalence relation
will be written $c_1(E,\U)$. The compatibility condition is
satisfied in for example equivalence relations given by the leaves
of a foliation on a manifold, and is the standard example of the
local equivalence relation defined by a foliation. \hfill
$\Box$\end{example}

We now consider similar questions for topological groupoids.

Of course if $G$ is a wide  subgroupoid of $Q$, then so also is
$C_1(G)$ and then $loc(C_1(G))$ is a local subgroupoid of $Q$.

Suppose $Q$ is {\it star path connected}, that is $Q= C_1(Q)$. Let
$X=\ob(Q)$ and let $U$ be a subset of $X$. In general $Q|U$ need
not be star path connected, as we show below. Further, while
$C_1(Q|U)\subs C_1(Q)|U$, in general we do not have equality here.
Such a condition is needed locally to obtain the local subgroupoid
$c_1(Q,\U)$ defined below.

\begin{defn} An open cover $\U=\{ U_i: i \in I \}$ of $X$
is said to be {\em path compatible } with a topological groupoid
$Q$ on $X$ if for all $i,j \in I, x \in U_i \cap U_j$ there is an
open set $W$ such that $ x \in W \subs U_i \cap U_j$ and
$$C_1(Q|U_i)|W = C_1(Q|U_j)|W.$$
In this case, the local subgroupoid  $c_1(Q,\U)$ is defined to
have value $[U_i,C_1(Q|U_i)]_x$ at $x \in U_i$.
\end{defn}
The next proposition gives useful sufficient conditions for
$c_1(Q,\U)$ to be defined.
\begin{prop} Let $Q$ be a  topological groupoid  on $X$ and suppose
there is an open cover $\mathcal{ U}=\{U_i\colon i\in I\}$ of $X$
such that for all $i,j \in I$ and $x \in U_i \cap U_j$ there is
an open set $W_x$ such that $ x \in W_x \subs U_i \cap U_j$ and
there are groupoid retractions $r_{i,W_x}: Q|U_i \to Q| W_x,
r_{j,W_x}:Q|U_j \to Q|W_x$ over retractions $ U_i \to W_x, U_j
\to W_x$. Then a local subgroupoid $c_1(Q,\U)$ is well defined by
for $i \in I, x \in U_i, x \mapsto [U_i,C_1(Q|U_i)]_x$.
\end{prop}
\begin{pf}
The retractions ensure the compatibility condition, since if $x,y
\in W$ and  if $\lambda$ is a path in $\St _{Q|U_i}$ joining
$1_x$ to the element  $g:x \to y$ of $Q|W$, then $r_{i,W}\lambda
$ is a path in $\St _{Q|W}$ joining $1_x$ to $g$. So
$C_1(Q|U_i)|W= C_1(Q|W)$, and similarly for $j$.
\end{pf}

Let $Q$ be a topological groupoid on $X$. Then $Q$ is called {\it
locally trivial} if for all $x \in X$ there is an open set $U$
containing $x$ and a  section $s\colon U\to \St _G x$ of $\beta$.
Thus $\beta s=1_{U}$ and  for each $y\in U$, $\alpha(s(y))=x$,
i.e. $s(y): x \to y$ in $Q$. We recall the following standard
result (see for example \cite{Ma}).

\begin{prop}\label{Mac} Let $ Q$ be a topological groupoid on $X$ and $U$ be
an open subset of $X$. If $s\colon U\rightarrow \St _Q x$ is a
continuous section of $\beta$ for some $x\in U$, then the
topological groupoid $Q|U$ is topologically isomorphic to the
product  groupoid $Q(x)\times (U\times U)$, and if $x \in W \subs
U$, then any retraction $U \to W$ is covered by a retraction $Q|U
\to Q|W$.
\end{prop}
\begin{pf}
Remark that the groupoid multiplication on $Q(x)\times (U\times
U)$ is defined by \[ (g,(y,z))(h,(z,w))=(gh,(y,w)).\] Define
\[\phi\colon Q|U\rightarrow G(x)\times (U\times U),
 g\mapsto (s(y)gs(z)^{-1},(y,z))\]
where $y=\alpha(g)$ and $z=\beta(g)$. Since $s$ is continuous,
$\phi$ is clearly an isomorphism of topological groupoids.

The last part follows easily.
\end{pf}

We have emphasised these results, despite their simple proofs,
because they have useful applications for example to manifold and
bundle theory.

If $s$ is a local subgroupoid of $Q$ defined by an atlas
$\U=\{(U_i,H_i):i \in I\}$ and $U$ is an open subset of $X$ then
$s|U$ is the local subgroupoid of $Q|U$ defined by the atlas
$\U\cap U=\{(U_i\cap U,H_i|(U_i\cap U)) :i \in I\}$. It is easy
to verify this is an atlas, and as a section $s|U$ is just the
restriction of $s$ to the open subset $U$.

Suppose now that we have the local subgroupoid $c_1(Q,\U)$
defined by the open cover $\U$, and $U$ is an open subset of $X$.
We will later need a result which  follows easily from
compatibility :
\begin{propdef} \label{pathlocal}
The equality
$$c_1(Q,\U)|U= c_1(Q|U,\U \cap U)$$
holds if for any $i,j \in I$ and $x \in U_i \cap U_j \cap U$ there
is an open set $W$ such that $x \in W \subs U_i \cap U_j \cap U$
and $C_1(Q|U_i)|W= C_1(Q|U_j \cap U)|W$. If this condition holds
for all open sets $U$ of $X$, then we say that the cover $\U$ is
{\em path local} for $c_1(Q,\U)$.\hfill $\Box$
\end{propdef}

\begin{rem}   There is a variation of the local subgroupoid
$c_1(Q,\U)$ in which the  paths in $Q$ which are used are
controlled, for example to belong to a given class, or to derive
from the paths in $X$ in a specified way. We give an example of
this in the next section.
 \hfill $\Box$
\end{rem}

\section{Path connections}

The purpose of this section is to give new examples of local
subgroupoids with a possibility of working towards relating the
concepts of holonomy  in foliation theory and in bundle theory.

Let $\Lambda (X)$ denote the path space of a topological space
$X$. Let $Q$ be a topological groupoid over $X$. A {\em path
connection} \cite{Br-S,Ma,Vi} $\Gamma$ in $Q$ is a continuous map
\[\Gamma \colon \Lambda(X)\rightarrow \Lambda (Q), \lambda\mapsto
\Gamma(\lambda)\] satisfying the following conditions
\begin{enumerate}[(i)]
  \item $\alpha(\Gamma(\lambda)(t))=\lambda(0)$  and
$\beta(\Gamma(\lambda)(t))=\lambda(t)$,  $t\in [0,1]$
  \item the transport  condition: If \[ \psi \colon [0,1]\rightarrow [t_0,t_1]\subseteq [0,1]\]
is a homeomorphism, then
\[\Gamma(\lambda)\circ\psi=\Gamma(\lambda)(\psi(0))\circ \Gamma(\lambda\psi).\]
\end{enumerate}
The second condition means
\[\Gamma(\lambda)(\psi(t))=\Gamma(\lambda)(\psi(0))\circ
\Gamma(\lambda\psi)(t) \] for $t\in [0,1]$,

By taking the homeomorphism $\psi$ to be the identity map
$\psi\colon [0,1]\rightarrow [0,1]$ it follows from the condition
(ii) that $\Gamma(\lambda)(0)=1_{\lambda(0)}.$ Let  $\lambda, \mu
\in \Lambda(X)$ and $\lambda(1)=\mu(0)$, that is the composition
$\lambda+\mu$ is defined, then we have
$\lambda=(\lambda+\mu)\circ\psi_0$ and
$\mu=(\lambda+\mu)\circ\psi_1$ where $\psi_0(t)=\frac{1}{2}t$ and
$\psi_1(t)=\frac{1}{2}t+\frac{1}{2}$. Moreover applying (ii) to
the path $\lambda+\mu$ and $\psi_0$ and then applying to $\lambda
+ \mu$ and $\psi_1$ we obtain
\begin{align*}
\Gamma(\lambda+\mu)(t)&=\left\{\begin{array}{ll}
\Gamma(\lambda)(2t)&\textrm{$0\leqslant t\leqslant \frac{1}{2}$}\\
\Gamma(\lambda)(1)\circ\Gamma(\mu)(2t-1)&
\textrm{$\frac{1}{2}\leqslant
  t\leqslant 1$}.
\end{array}\right. \\
\intertext{In particular}
\Gamma(\lambda+\mu)(1)&=\Gamma(\lambda)(1)\circ\Gamma(\mu)(1).
\end{align*}
Let $Q$ be a topological groupoid on $X$ with a continuous path
connection $\Gamma \colon \Lambda(X)\rightarrow \Lambda (Q).$ Let
$C_{\Gamma}(Q)$ be the set of all $g\in Q$ such that if
$\alpha(g)=x$ then there is a path $\lambda$ in $X$ such that
$\Gamma(\lambda)$ joins $g$ to the identity $1_{x}$ at $x$, that
is, $\Gamma(\lambda)(0)=1_x$ and $\Gamma(\lambda)(1)=g$. Then we
prove the following proposition.
\begin{prop}
$C_{\Gamma}(Q)$ is a wide subgroupoid of $Q$.
\end{prop}
\begin{pf} Let $g,h\in C_{\Gamma}(Q)$ such that $gh$ is defined in $Q$. Then there
are paths $\lambda$ and $\mu$ joining $g$ to $1_{\alpha(g)}$ and $h$ to
$1_{\alpha(h)}$ respectively. Here note that $\lambda(0)=\alpha(g)$,
 $\lambda(1)=\beta(g)$ and  $\mu(0)=\alpha(h)$,
 $\mu(1)=\beta(h)$. So the composition $\lambda+\mu$ of the paths is defined
and $\Gamma(\lambda+\mu)(0)=\Gamma(\lambda)(0)=1_{\alpha(g)}$ and
$\Gamma(\lambda+\mu)(1)=\Gamma(\lambda)(1)\circ\Gamma(\mu)(1)=gh$.
So $gh\in C_{\Gamma}(Q)$. That means  $C_{\Gamma}(Q)$ is closed
under the groupoid composition.

If $g\in C_{\Gamma}(Q)$ with $\alpha(g)=x$ then there is a path $\lambda$ in
 $X$ such
that $\Gamma(\lambda)(0)=1_x$ and $\Gamma(\lambda)(1)=g$. Define a path
$\bar{\lambda}$ in $X$ such that $\bar{\lambda}(t)=\lambda(1-t)$. Then
$\bar{\lambda}(t)=(\lambda\psi)(t)$ with $\psi(t)=1-t$.
By the transport law we have
 $ \Gamma(\lambda)(\psi(t))=\Gamma(\lambda)(1)\circ\Gamma(\bar{\lambda})(t)$
where $\Gamma(\lambda)(1)=g$. So we have
\[ \Gamma(\bar{\lambda})(0)=g^{-1}\Gamma(\lambda)(1)=g^{-1}g=1_y  \]
and
\[ \Gamma(\bar{\lambda})(1)=g^{-1}\circ(\Gamma(\lambda)(0)=g^{-1}\circ 1_x=g^{-1}
\]
So $g^{-1}\in C_{\Gamma}(Q)$. Hence $C_{\Gamma}(Q)$ is a wide subgroupoid of $Q$.
\end{pf}

We also need an analogue of Proposition \ref{generate}.

\begin{prop} \label{conn-gen}
If $\Gamma$ is a path connection on the topological groupoid $Q$
and $\U$ is an open cover of $X$, then $C_\Gamma(Q)$ is generated
by the family $C_\Gamma(Q|U)$ for all $U \in \U$.
\end{prop}
\begin{pf}
If $g \in C_\Gamma(Q)$ is joined to $1_x$ by the path
$\Gamma(\mu)$, then we can write $\mu = \mu_1+\cdots+\mu_n$ where
each $\mu_r$ lies in some set $U_r$ of $\U$. Let $g_r =
-\Gamma(\mu_r)(0)+ \Gamma(\mu_r)(1)$. Then $g_r \in
C_\Gamma(Q|U_r)$ and $g = g_1+\cdots+g_n$.
\end{pf}

 If $Q$ is  a topological groupoid on $X$ with a path connection
$\Gamma\colon \Lambda(X)\rightarrow \Lambda(Q)$ then of course
$loc(C_{\Gamma}(Q))$ is a local subgroupoid. However we would like
an analogue of $c_1(Q,\U)$ and this needs extra conditions. In
fact the existence of a smooth path connection for the groupoid
$\it Sym(p)$ of a principal bundle $p\colon E\rightarrow B$
relies on the existence of an infinitesimal connection (see
\cite{Kob-Nom}, \cite{Ma}) which itself requires extra structure
on the space involved.

We give some conditions which are sufficient for
$c_{_\Gamma}(Q,\U)$ to be well defined.

We suppose given an open cover $\mathcal U=\{U_i\colon i\in I\}$
for $X$ and for each $i\in I$ a collection $geod(U_i)$ of paths
in $U_i$ --  an element $\lambda\in geod(U_i)$ with
$\lambda(0)=x$, $\lambda(1)=y$ is called a ``geodesic path'' from
$x$ to $y$. We suppose \begin{enumerate}[(i)] \item  if $x,y\in
U_i$, then there is a unique geodesic path $geod_i(x,y)$ from $x$
to $y$. \item  if $x,y\in U_i\cap U_j$ then
$geod_i(x,y)=geod_j(x,y)$. \end{enumerate}

We also need the connection to be ``flat'' for this structure in
the sense that   if $\lambda\colon x\rightarrow y$ is any path in
$U_i$ then $ \Gamma(\lambda)(1)=\Gamma(geod_i(x,y))(1)$.

\begin{prop}Under the above atlas  assumptions, there is a local subgroupoid
$c_{_\Gamma}(Q,\U)$ defined by
\[c_{_\Gamma}(Q,\U)(x)=[U_i,C_{\Gamma}(Q|U_i)]_x.\]
\end{prop}
\newpage
\begin{pf} We have to prove that if $x\in U_i\cap U_j$ then
\[ [U_i,C_{\Gamma}(Q| U_i)]_x= [U_j,C_{\Gamma}(Q| U_j)]_x.\]
This means there is an open neighbourhood $W$ of $x$ in $ U_i\cap
U_j$ such that
\[  C_{\Gamma}(Q| U_i)| W= C_{\Gamma}(Q| U_j)| W.\]
Let $W$ be an open neighbourhood of $x$ in $U_i\cap U_j$. Let
$g\in C_{\Gamma}(Q| U_i)| W$ with $\alpha(g)=x$ and $\beta(g)=y$.
So there is a path $\lambda\colon x\rightarrow y$ in $U_i$ such
that $\Gamma(\lambda)(1)=g$. Let $\lambda_i\colon x\rightarrow y$
be the geodesic path. So $\Gamma(\lambda_i)(1)=g$, by the flat
condition,  and so $g\in C_{\Gamma}(Q| U_j)| W$. Hence $
C_{\Gamma}(Q| U_i)| W\leq C_{\Gamma}(Q| U_j)| W$. Since the
converse proof is similar we have $C_{\Gamma}(Q| U_i)|
W=C_{\Gamma}(Q| U_j)| W$.
\end{pf}

Suppose now that we have the local subgroupoid $c_{_\Gamma}(Q,\U)$
defined by the open cover $\U$, and $U$ is an open subset of $X$.
We will later need a result which  follows easily from
compatibility :
\begin{propdef} \label{g-pathlocal}
The equality
$$c_{_\Gamma}(Q,\U)|U= c_{_\Gamma}(Q|U,\U \cap U)$$
holds if for any $i,j \in I$ and $x \in U_i \cap U_j \cap U$ there
is an open set $W$ such that $x \in W \subs U_i \cap U_j \cap U$
and $C_{\Gamma}(Q|U_i)|W= C_{\Gamma}(Q|U_j \cap U)|W$. If this
condition holds for all open sets $U$ of $X$, then we say that
the cover $\U$ is {\em  $\Gamma$ path local} for
$c_{_\Gamma}(Q,\U)$.\hfill $\Box$
\end{propdef}

\section{Partial orders and induced morphisms}
We first establish  some elementary but essential basic theory.

The set $L_Q(X)$ of wide subgroupoids of $Q$ is a poset under
inclusion. We write $\leq$ for this partial order. This poset has
a top element namely $Q$ and a bottom element namely the discrete
groupoid $X$.

Let ${\Loc}(Q)$ be the set of local subgroupoids of $Q$. Let $x\in
X$. We define a partial order on the stalks ${p_Q}^{-1}(x)=(\L
_Q)_x$ by $[U', H']_x\leq [U, H]_x$ if there is an  open
neighbourhood $W$ of $x$ such that $W\subseteq U\cap U'$ and
$H'|W$ is a subgroupoid of $H|W$. Clearly this partial order is
well defined. Its bottom element is of the form $[U,H]_x$ where
$H$ is discrete, and its top element is of the form $[U,Q|U]_x$.
This partial order induces a partial order on ${\bf Loc}(Q)$ by
$s\leq t$ if and only if $s(x)\leq t(x)$ for all $x\in X$.

The major purpose of the next topic is to relate local
subgroupoids of $Q$ and local equivalence relations on $\ob(Q)$.
This seems an area requiring much more development, and we hope
will be the start of new ways of  relating bundle and foliation
theory.

Suppose given two groupoids $Q, H$ and a groupoid morphism $\phi
\colon Q\rightarrow H$ on $X$, which is the identity on objects.
Then we obtain morphisms of presheaves $\phi^*: L_Q \to  L_H $,
$\phi_*:L_H\to L_Q$ as follows.

Let $U$ be open in $X$. Then $\phi_* \colon L_Q(U)\rightarrow
L_H(U)$ is given by $\phi(K)$ is the image of $K\in L_Q(U)$ by
$\phi$. Here $K$ is a wide subgroupoid of $Q|U$, and so  its image
$\phi(K)$ is a subgroupoid of $H|U$, since $Ob(\phi)$ is
injective and is clearly wide.

Further $\phi^* \colon L_H(U)\rightarrow L_Q(U)$ is given by
$\phi^*(K')=\phi^{-1}(K')$, for $K'\in L_H(U)$.

Hence we get induced morphism of sheaves $\phi_*:\L _Q\to \L _H,$
$\phi^*:\L _H\to \L _Q$.

In particular, we get for a groupoid $Q$ an `anchor' morphism of
groupoids $A: Q\to X\times X$ and so sheaf morphisms
\[ A_*: \L _Q\to \L _{X\times X}, \ \ \
A^* : \L _{X\times X}\to \L _Q. \] Hence a local subgroupoid $s$
of $Q$ yields a local equivalence relation $A_*(s)$ on $X$, and a
local equivalence relation $r$ on $X$ yields a local subgroupoid
$A^*(r)$ of $Q$. This gives further examples of local
subgroupoids.

Clearly also $\phi_*$, $\phi^*$ are order preserving on stalks for
any morphism $\phi:Q\to H$ of groupoids over $X$. Hence they
induce morphism of posets
\[ \phi_*:{\bf Loc}(Q)\to {\bf Loc}(H), \ \ \
\phi^*: {\bf Loc}(H)\to {\bf Loc}(Q).  \]
Further, $s\leq \phi^*r$ if and only if $\phi_*s\leq r$.
This can be expressed by saying that $\phi_*$ is left adjoint to
$\phi^*$.

%----------------------------------------------------------------------------------
\section{Coherence for wide subgroupoids and local subgroupoids}

 We now fix a groupoid  $Q$  on $X$, so that $L_Q(X)$ is
the set of wide subgroupoids of $Q$, with its inclusion partial
order, which we shall write $\leq$.

Clearly $loc_Q$ as defined in Definition \ref{loc} gives a poset morphism
\[loc_Q\colon L_Q(X)\rightarrow \Loc(Q).\]
\begin{defn}
Let $s$ be a local subgroupoid of  $Q$. Then $glob(s)$ is the wide
subgroupoid of $Q$ which is the intersection  of all wide
subgroupoids $H$ of $Q$ such that $s\leq loc(H)$. \hfill $\Box$
\end{defn}
We think of $glob(s)$ as an approximation to $s$ by a global
subgroupoid.
\begin{prop}\label{sim}
\begin{enumerate}[\rm (i)]
\item  $loc$ and $glob$ are morphisms of posets.
\item  For any wide subgroupoid $H$ of $Q$, $glob(loc(H))\leq H$.\hfill $\Box$
\end{enumerate}
\end{prop}
The proofs are clear.

However, $s\leq loc(glob(s))$ need not hold. Rosenthal in
\cite{Ro1} gives the example of the local equivalence relation
$s=loc(E)$ where $E$ is the equivalence relation $aEb$ if and
only if $a =\pm b$. Here is a similar example.
\begin{example}
Let $Q$ be a groupoid on $\Reno$ such that all $x \in \Reno$ with
$x \ne 0$ there is a neighbourhood $U$ of $x$ such that $Q|U$ is a
bundle of groups, while no such neighbourhood of $0$ exists. Let
$s=loc(Q)$. Then $H=glob(s)$ on $\Reno \setminus \{0\}$ coincides
with $Q$ on this set, and in fact $H$ is the bundle of groups
$Q(x)$ for all $ x \in \Reno$. It follows that $s(0)> loc(H)(0)$.
\hfill $\Box$
\end{example}

We therefore adapt from \cite{Ro1,Ro2} some notions of coherence.

\begin{defn}
Let $s$ be a local subgroupoid of $Q$ on $X$.
\begin{enumerate}[(i)]
  \item   $s$ is called {\bf coherent} if  $s\leq loc(glob(s))$.
\item $s$ is called {\bf globally coherent} if $s = loc(glob(s))$.
\item  $s$ is called {\bf totally coherent} if for every open set $U$ of $X$, $s|U$
is coherent.
\end{enumerate}
\end{defn}
\begin{defn}
Let $H\in L_Q(X)$, so that $H$ is  a wide subgroupoid of $Q$.
\begin{enumerate}[(i)]
  \item  $H$ is called {\bf locally coherent} if $loc(H)$ is coherent.
  \item  $H$ is called {\bf coherent} if $H = glob(loc(H))$.
  \hfill $\Box$
\end{enumerate}
\end{defn}

\begin{example}
Let $Q$ be a groupoid on the discrete space $X$. Then
$glob(loc(Q))= Inn(Q)$, the groupoid of vertex groups of $Q$. Thus
in general, $Q$ is not coherent. \hfill $\Box$
\end{example}
At another extreme we have:
\begin{prop}
Let $Q$ be a bundle of groups. Then any local subgroupoid of $Q$
is globally coherent, and any wide subgroupoid of $Q$ is coherent.
\end{prop}
\begin{pf}
Let $s$ be a local subgroupoid of $Q$ and let $\{(U_i, H_i) :i\in
I\}$ be an atlas for $s$. Then if $x\in U_i$, we have $s(x)=[U_i,
H_i]_x$. Let $H(x)=H_i(x)$. If $x\in U_i\cap U_j$, there is a
neighbourhood $W$ of $x$ such that $W\subseteq U_i\cap U_j$ and
$H_i|W = H_j|W$, and hence $H_i(x) = H_j(x)$. Thus $H$ is
independent of the choices. Also $H|U_i = H_i$. Hence $loc(H)(x) =
[U_i, H_i]_x$, and so $loc(H)=s$, $H= glob(s)$.
\end{pf}

Coherence of $s$ says that in passing between local and global
information nothing is lost due to collapsing. Notice also that
these definitions depend on the groupoid $Q$.
\begin{prop}
$loc$ and $glob$ induce  morphisms of posets from the locally
coherent subgroupoids of $Q$ to the coherent local subgroupoids of
$Q$, and on these posets $glob$ is left adjoint to $loc$. Further,
$glob$ and $loc$ give inverse  isomorphisms between the posets of
coherent  subgroupoids and of  globally coherent local
subgroupoids.
\end{prop}
\begin{pf}
Let $H$ be a locally coherent subgroupoid of $Q$  and let $s$ be a
coherent local subgroupoid of $Q$. By the definition of locally
coherent subgroupoid $loc(H)$ is a coherent local subgroupoid of
$Q$.

Conversely, let $K=glob(s)$. Since $s$ is coherent, $s\leq
loc(K)$. Since $glob$ is a poset morphism, $glob(s)\leq
glob(loc(K))$, i.e., $K\leq glob(loc(K))$. Since $loc$ is also a
poset morphism,
\[ loc(K)\leq loc(glob(loc(K))). \] So
$loc(K)$ is a coherent local subgroupoid, and $K$ is locally
coherent.

The adjointness relation is that $ glob(s)\leq H
\Longleftrightarrow s \leq loc(H).$ The implication $\Leftarrow$
follows from the fact that for all $H$ we have $ glob(loc(H)) \leq
H$. The implication $\Rightarrow$ follows from the coherence of
$s$.

The final statement is obvious.
\end{pf}

Note in particular that coherence of $H$ implies local coherence
of $H$.

\begin{prop}\label{coh1}
Let $Q$ be a topological groupoid on $X$ and  $G$ a star path
connected wide subgroupoid of $Q$. Then $G$ is coherent and
$loc(G)$ is globally coherent.
\end{prop}
\begin{pf} We prove that $glob(loc(G))=G$. By Proposition \ref{sim}
 $glob(loc(G))\leq G$. To prove that $ G\leq glob(loc(G))$ let $H$ be
a wide subgroupoid of $Q$ such that $loc(G)\leq loc(H)$. Then for
$x\in X$,
\[[X,G]_x\leq [X,H]_x\]
and so for some open neighbourhood $U_x$ of $x$, $G|U_x\leq
H|U_x$. These sets $U_x$ form a cover $\U$ of $X$. By Proposition
\ref{generate}, $G$ is generated by the $G|U$ for $U \in \U$.  It
follows that $G \leq H$.
\end{pf}

\begin{cor}
If $X$ is a topological space then its fundamental groupoid
$\pi_1X$ is coherent and globally coherent.
\end{cor}
\begin{cor} Let $Q$ be a topological groupoid on $X$. Then the star
identity component $C_1(Q)$ is coherent.
\end{cor}

\begin{cor} Let $Q$ be a topological groupoid with a path connection $\Gamma\colon
\Lambda(X)\rightarrow \Lambda(Q)$. Then the wide subgroupoid
$C_{\Gamma}(Q)$ is coherent.
\end{cor}

\begin{prop}
Let $Q$ be a topological groupoid on $X$. Suppose that the local
subgroupoid $c_1(Q,\U)$ of $Q$ is well defined by an open cover
$\U$. Then
\begin{enumerate}[\rm (i)]
  \item $glob(c_1(Q,\U))=C_1(Q)$.
  \item $c_1(Q,\U)$ is coherent.
  \item If $\U$ is path local for $Q$, then $c_1(Q,\U) $ is  totally coherent.
\end{enumerate}
\end{prop}
\begin{pf}
\noindent (i) Certainly $c_1(Q,\U) \le loc(C_1(Q))$ since for all
open $U$ in $X$, $C_1(Q|U) \leq C_1(Q)|U$ and so
$[U_i,C_1(Q|U_i)]_x \le [X,C_1(Q)]_x$ for all $x \in U$.

Now suppose $H$ is a wide subgroupoid of $Q$ and $c_1(Q,\U) \le
loc(H)$. We have to prove $C_1(Q) \le H$.

Let  $i \in I$ and $x \in U_i$. We have   $ [U_i,C_1(Q|U_i)]_x \le
[X, H]_x$.  Hence  there is an open neighbourhood  $W_x$ of $x$
contained in $U_i$ and such that $$C_1(Q|U_i)|W_x \leq H.$$ By
Proposition \ref{generate}, $C_1(Q|U_i)$ is generated by the
$C_1(Q|U_i)|W_x$ for all $ x \in U_i$ and by the same Proposition,
$C_1(Q)$ is generated by the $C_1(Q|U_i)$ for all $i \in I$. Hence
$C_1(Q)\leq H$.

Coherence of  $c_1(Q,\U) $ follows from  (i)  and $C_1(Q|U_i) \leq
C_1(Q)|U_i$. Total coherence in the path local case  follows by
applying (ii)  to $Q|U$, using Proposition and Definition
\ref{pathlocal}.
\end{pf}

\begin{prop} Let $Q$ be a topological groupoid on $X$ such that
the local subgroupoid $c_{_\Gamma}(Q,\U)$ is well defined by the
open cover $\U$. Then:
\begin{enumerate}[\rm (i)]
  \item $glob(c_{_\Gamma}(Q,\U))=C_{\Gamma}(Q)$.
  \item $c_{_\Gamma}(Q,\U)$ is coherent.
  \item If $\U$ is $\Gamma$ path local for $Q$, then $c_{_\Gamma}(Q,\U)$
  is totally coherent.
\end{enumerate}
\end{prop}\label{coh3}
\begin{pf} (i) Note that $c_{_\Gamma}(Q,\U)\leq loc(C_{\Gamma}(Q))$ since for all $U$
  in $X$, $C_{\Gamma}(Q|U)\leq C_{\Gamma}(Q)|U$. So
  $glob(c_{\Gamma}(Q))\leq C_{\Gamma}(Q)$.
To prove that  $ C_{\Gamma}(Q)\leq glob(c_{_\Gamma}(Q,\U))$
suppose that $H$ is a wide subgroupoid of $Q$ such that
$c_{_\Gamma}(Q,\U)\leq loc(H)$. We have to prove that
$C_{\Gamma}(Q)\leq H$. If  $U \in \U$ and $x\in U$ then
\[[U,C_{\Gamma}(Q| U)]_x\leq [X,H]_x.\]
Hence $U$ has a covering by open sets $W_x$ such that $
[U,C_{\Gamma}(Q|U)|W_x \leq  H$. By Proposition \ref{generate},
$C_{\Gamma}(Q|U)$ is generated by the groupoids $
C_{\Gamma}(Q|U)|W_x$ and by Proposition \ref{conn-gen}
$C_{\Gamma}(Q)$ is generated by the $C_{\Gamma}(Q|U)$ for $U \in
\U$. Hence $C_{\Gamma}(Q)\leq H$.

The proofs of (ii), (iii) are analogous to those in the previous
proposition.

\end{pf}

\section{Coherence and atlases }

We lead up to conditions for an atlas for a local subgroupoid $s$
to be {\em globally adapted} to $s$. This notion is important for
considerations of holonomy (see \cite{Br-Ic}), and the
applications will be developed elswhere.

The next proposition gives an alternative description of $glob$.

Let ${\cal U}_s = \{(U_i, H_i):i\in I\}$ be an atlas for the local
subgroupoid $s$. Then $glob({\cal U}_s)$ is defined to be the
subgroupoid of $Q$ generated by all the $H_i, i \in I$.

An atlas ${\cal V}_{s}=\{(V_j, {s}_j): j\in J\}$ for $s$ is
said to refine ${\cal U}_s$ if for each index $j\in J$
there exists an index $i(j)\in I$ such that $V_j\subseteq U_{i(j)}$
and $s_{i(j)}|V_j = s_j$.

\begin{prop} \label{refin2}
Let  $s$  be a local subgroupoid of $Q$ given by the atlas ${\cal
U}_s = \{(U_i, H_i):i\in I\}$. Then $glob(s)$ is the intersection
of the subgroupoids $glob({\cal V}_s) $ of $Q$ for all refinements
${\cal V}_s$ of ${\cal U}_s$.
\end{prop}
\begin{pf}
Let $K$ be the intersection given in the proposition.

Let $Q$ be a subgroupoid of $Q$ on $X$ such that $s\leq loc(Q)$.
Then for all $x \in X$ there is a neighbourhood $V$ of $x$ and
$i_x \in I$ such that $x \in U_{i_x}$ and $ H_{i_x}|V_x \cap
U_{i_x} \leq Q$. Then ${\cal W}= \{(V_x \cap U_{i_x},H_{i_x}|V_x
\cap U_{i_x}) : x \in X\} $ refines ${\cal U}_s$ and $glob({\cal
W}) \leq Q$. Hence $K \leq Q$, and so $ K \leq glob(s) $.

Conversely, let ${\cal V}_s = \{ (V_j,H'_j) : j \in J \}$ be an
atlas for $s$ which refines ${\cal U}_s$. Then for each $ j \in J$
there is an $i(j) \in I$ such that $ V_j \subseteq
U_{i(j)},H'_j=H_{i(j)}|V_j$. Then $s \leq loc(glob({\cal V}_s))$. Hence
$glob(s) \leq glob({\cal V}_s)$ and so $glob(s) \leq K$.
\end{pf}

\begin{cor} \label{cohgen}
A wide  subgroupoid $H$ of $Q$ is  coherent if and only if
for every open cover ${\cal V} $ of $X$,  $H$   is
generated by the subgroupoids $ H|V , V \in {\cal V} $.
\end{cor}
\begin{pf}
Note that $\{(X,H)\}$ is an atlas for $loc(H)$, which is refined by
${\cal V}_H= \{(V,H|V): V \in \cal V\}$ for any open cover $\cal V$.

Suppose the latter condition holds. Then Proposition \ref{refin2}
implies that $H = glob(loc(H))$, i.e. $H$ is coherent. The converse
holds since $glob({\cal V}_H) \leq H$.
\end{pf}

Let $U$ be an open subset of $X$. Then we have notions of local
subgroupoids of $Q|U$ and also of the restriction $s|U$ of a local
subgroupoid $s$ of $Q$. Clearly if $H$ is a wide subgroupoid of
$Q$ then $loc(H|U)= (loc(H))|U$.
\begin{prop} \label{subres}
Let $s$ be a local subgroupoid of $Q$ and let $U$ be open in $X$.
Then $glob(s|U)\leq glob(s)|U$.
\end{prop}
\begin{pf}
Let $H$ be a wide subgroupoid of $Q$ such that $s \leq loc(H)$.
Then $s|U \leq loc(H|U)$.  So $glob(s|U) \leq H|U$. The result
follows.
\end{pf}

\begin{prop}\label{pr1} Let $s$ be a local subgroupoid of $Q$.
Then
\begin{enumerate}[\hspace{-1em} \rm (i)]
  \item If $s$ is globally  coherent, $U$ is open in $X$, and $s|U$ is
coherent, then $s|U$ is globally coherent.
  \item If there is an open cover $\cal V$ of $X$ such that $s|V $ is
 coherent for all $V \in \cal V$, then $s$ is coherent.
  \item  If $s$ is globally coherent then for any open cover $\cal V$
of $X$,  $glob(s)$  is generated by the groupoids $glob(s)|V$ for
all $V \in \cal V$.
\item If there is open cover ${\cal V}$ of $X$ such that
$s|V$ is globally and totally coherent for $V\in {\cal V}$, then
$s$ is totally coherent.
\end{enumerate}
\end{prop}
\begin{pf}
i) We are given $s = loc(glob(s))$. By Proposition \ref{subres}
\[ loc(glob(s|U))\leq loc(glob(s)|U) = loc(glob(s))|U =
s|U.\]
Since $s|U$ is coherent, we have $s|U\leq loc(glob(s|U))$.
So $s|U = loc(glob(s|U))$, i.e. $s|U$ is globally coherent.

ii) We have
$$s|V \leq loc(glob(s|V)) \leq loc(glob(s)|V) \leq
(loc(glob(s)))|V .$$
Since this holds for all $V$ of an open cover, we have $s \leq
loc(glob(s)).$

iii) This follows from Corollary \ref{cohgen}.

iv) Let $U$ be open in $X$. Let $V\in {\cal V}$.
Since $s|V$ is globally and totally coherent,
then $s|V\cap U$ is globally coherent. Hence  by (ii) $s|U$ is coherent,
since the $V\cap U, V\in {\cal V}$, cover $U$.
\end{pf}

\begin{prop}
Let ${\cal U}_s=\{ (U_i,H_i):i \in I\}$ be an atlas for the local
subgroupoid $s$. Then:
\begin{enumerate}[\hspace{-0.5em}\rm (i)]
  \item  $s|U_i = loc(H_i)$ for all $ i \in I$;

  \item  $loc(glob(s|U_i)) \leq s|U_i$ for all $ i \in I$;

  \item if $s|U_i$ is coherent for all $ i \in I$ then $s$  is globally coherent;
\item if $s|U_i$ is coherent for all $i \in I$ then $glob(s)=glob({\cal U}_s)$.
\end{enumerate}

\end{prop}
\begin{pf}
i) This is clear.

ii) We have  by Proposition \ref{sim}
$$ loc(glob(s|U_i)) = loc(glob(loc(H_i))) \leq loc(H_i) = s|U_i. $$

iii) This is immediate from the definition of coherence and ii).

iv) Let $H=glob({\cal U}_s)$, i.e. $H$ is the subgroupoid of $Q$
generated by the $H_i, i \in I$. Then $glob(s) \leq H$. Let $K$ be
a wide subgroupoid of $Q$ such that $s \leq loc(K)$. Then for all
$i \in I$ and $ x \in U_i$ there is a neighbourhood $V_x^i$ of $x$
such that $V_x^i \subseteq U_i$ and $H_i|V_x^i \leq  K|V_x^i$. By
global coherence of $s|U_i$ and Proposition \ref{pr1}(i,ii), $H_i$
is generated by the $H_i|V_x^i$ for all $x \in U_i$. Hence $H_i
\leq K|U_i\leq K$. Hence $H \leq K$. Hence $H \leq glob(s)$.
\end{pf}

\begin{defn}\label{def1}
Let $s$ be a local subgroupoid of the groupoid $Q$ on $X$. An
atlas ${\cal U}_s$ for $s$ is called {\it  globally adapted} if $
glob(s)= glob({\cal U}_s)$.
\end{defn}

\begin{rem}This is a variation on the notion of an
$r$-adaptable family defined by Rosenthal in \cite[Definition
4.4]{Ro2} for the case of a local equivalence relation $r$. He
also imposes a connectivity condition on the local equivalence
classes.
\end{rem}

The construction of the holonomy groupoid of a local subgroupoid
requires a globally adapted atlas (see Theorems 3.7 and 3.8
\cite{Br-Ic}). The following proposition is very useful to this
end.
\begin{prop}
Let $s$ be a totally coherent local subgroupoid of the groupoid
$Q$ on $X$. Then any atlas for $s$ is    globally adapted.
\end{prop}
\begin{pf}
This is immediate from the previous Proposition, since total
coherence implies that each $s|U_i$ is coherent.
\end{pf}
\begin{cor}\label{glo-adapt} Any path local atlas $\U$ of the local subgroupoid $c_1(Q,\U)$
is globally adapted.
\end{cor}
\begin{cor} Any $\Gamma$ path local atlas of the local subgroupoid $c_{_\Gamma}(Q,\U)$
is globally adapted.
\end{cor}

{}

\end{document}